\documentclass{amsart}
\usepackage[latin1]{inputenc}
\usepackage{amsmath}
\usepackage{amsthm}
\usepackage{amsfonts}
\usepackage{amssymb}
\newcommand{\ep}{\hspace*{\fill}$\Box$}

\begin{document}
\title[$P(f)=Q(g)$]{Rational Decompositions of Complex Meromorphic Functions}
\author[A. Escassut and E. Mayerhofer]{Alain Escassut and Eberhard Mayerhofer}
\address{University of Vienna, Department of Mathematics, Nordbergstrasse 15, 1090 Vienna, Austria}
\email{eberhard.mayerhofer@univie.ac.at}
\begin{abstract} Let $h$ be a  complex meromorphic function. The 
problem of decomposing $h$ in two different ways
$P(f)$ and
$Q(g)$ with
$f,\ g$ two other  meromorphic functions and $P,\ Q$ polynomials,     was 
studied by C.-C. Yang,  P. Li and K. H. Ha. Here we
consider the problem when we replace the polynomials $P,\ Q$ by rational 
functions $F,\ G$. Let $\deg(F)$ be the maximum
degree of numerator and denominator of $F$. Assume  some zeros $c_1,\dots,c_k$ 
of $F'$ satisfy a pack of 5 conditions involving
particularly $G(d)\neq F(c_j)$ and $D(d)\neq 0$ for every zero $d$ of 
$C'-F(c_j)D'$, with $G=\frac{C}{D}$, ($j=1,\dots,k$).
First, we show  that  if
$f,\ g$ are entire functions such that $F(f)=G(g)$ then $k\deg(G)\leq \deg(F)$. 
Now, let $\gamma(D)$ be the number of distinct
zeros of the denominator of $G$ and assume that meromorphic functions $f,\ g$ 
satisfying $F(f)=G(g)$, then $k\deg(G)\leq
\deg(F)+k\gamma(D)$. When  zeros $c_1,\dots,c_k$ of $F'$ satisfy a stronger
condition, then we show that $k\deg(G)\leq
\deg(F)+k\min(\gamma(C), \gamma(D))$.\end{abstract}
\maketitle
\section{Introduction and results}
Let ${\mathcal A}(\mathbb{C})$ be the algebra of entire functions in $\mathbb{C}$ and ${\mathcal M}(\mathbb{C})$ be the field of meromorphic functions in
$\mathbb{C}$.  The decomposition of a function $h\in {\mathcal M}(\mathbb{C})$ was studied in \cite{2} and \cite{8} in order to find many cases when $h$
cannot admit two different decompositions of the form $P(f)$ and $Q(g)$, with $P,\ Q\in \mathbb{C}[x]$ and $f,\ g\in {\mathcal M}(\mathbb{C})$.  In \cite{1},  a similar study was made in  a  p-adic field $K$ of characteristic $p\geq 0$, concerning meromorphic
functions in all $K$, or meromorphic functions inside an ''open'' disk.  In \cite{6} the problem was considered in a  p-adic field again, when $P,\ Q$ lie not in $K[x]$ but in $K(x)$. Here we mean to consider this last  problem in $\mathbb{C}$: given $h\in {\mathcal M}(\mathbb{C})$,  can
we find conditions excluding that $h$ admits  two different decompositions of the form $F(f)$ and $G(g)$, with $F,\;G\in \mathbb{C}(x)$? The Nevanlinna Second Main Theorem provides us with the main tool to deal with the problem, in this case as in
the previous ones. We shall also use some basic algebraic consideration already made in [6]. 
\\{\bf Notation.} \quad Let $F=\frac{A}{B} \in \mathbb{C}(x)$   with $A,\;B  \in \mathbb{C}[x]$ and  $gcd(A,B)=1$. We put $\deg(F)=\max(\deg(A), \deg(B))$. 
Let $P\in \mathbb{C}[x]$. We denote by $\gamma(P)$ the number of distinct zeros of $P$. Given $F(x)= \frac{ A}{B}$  with  $A,\ B\in \mathbb{C}[x]$  and  $gcd(A,B)=1$, we put $\lambda(F)=\min(\gamma(A), \gamma(B))$.
\medskip\\{\bf Definition.}\quad Let $\mathbb{E}$ be an algebraically closed field, let  $F=\frac{A}{B},\ G=\frac {C}{D} \in \mathbb{E}(x)$ with $A, \ B ,\ C,\ D \in \mathbb{E}[x]$,  $gcd(A,B)=gcd(C,D)=1$. Let $c_1,\dots,c_k\in \mathbb{E}$ be zeros of $F'$. Then $(F,G)$ will be said {\it to satisfy Condition M for $c_1,\dots, c_k$} \ if:
\\\medskip  1) \ $C$ and $D$ are monic,
\\\medskip  2) \ $F(c_j)\neq F(c_l)\ \forall j\neq l$, 
\\\medskip  3) \ if $\deg(C)=\deg(D)$ we assume $F(c_j)\neq 1\ \forall j=1,\dots,k$,  
\\\medskip  4) \ for every $j=1,\dots,k$ we have $G(d)\neq F(c_j)$ and $D(d)\neq 0$ for every zero $d$ of $C'-F(c_j)D'$.\\ Further,  $(F,G)$ will be said {\it to satisfy Condition M' for $c_1,\dots, c_k$} \ if Conditions 1), 2), 3) are satisfied
and if Condition 4) is replaced by Condition 4'):\\\medskip 4')  \ for every $j=1,\dots,k$ we have $G(d)\neq F(c_j)$ and $C(d)D(d)\neq 0$ for every zero $d$ of $C'-F(c_j)D'$,   
\medskip\\ {\bf Theorem 1.}\quad {\it Let $F, \ G\in \mathbb{C}(x)$  satisfy Condition M for $c_1,\dots,c_k$  and let $p=\deg(F),\ q=\deg(G)$.     Assume that there exist two  functions $f,\ g\in {\mathcal A}(\mathbb{C})$ such that $F(f)=G(g)$. Then $kq\leq p$. }
\medskip\\  {\bf Theorem 2.}\quad {\it Let $F=\frac{A}{B},\ G=\frac{C}{D} \in \mathbb{C}(x)$ with $gcd(A,B)=gcd(C,D)=1$ and $C, D$
monic,  satisfy Condition M for $c_1,\dots,c_k$.  Let
$p=\deg(F),\ q=\deg(G)$.    
 Assume that there exist two  functions $f,\ g\in {\mathcal 
M}(\mathbb{C})$ such that $F(f)=G(g)$. Then $kq\leq p(1+(k\gamma(D))$. }
\medskip \\{\bf Theorem 3.}\quad {\it Let $F=\frac{A}{B},\ G=\frac{C}{D}\in \mathbb{C}(x)$ with $gcd(A,B)=gcd(C,D)=1$ and $C, D$
monic,  satisfy Condition M' for $c_1,\dots,c_k$.  Let
$p=\deg(F),\ q=\deg(G)$.    
 Assume that there exist two  functions $f,\ g\in {\mathcal 
M}(\mathbb{C})$ such that $F(f)=G(g)$. Then $kq\leq p(1+k\lambda(G))$. }
\medskip\\ {\bf Remarks.} \quad 1)  \ Theorems 2 and 3 are trivial when $C$ and $D$ have same degree and have no multiple zero, and
are more interesting when
$\min(\deg(C),\deg(D))$ is much smaller than $\max(\deg(C),\deg(D))$. 
\medskip\\ 2)\   Let $f(x)=\sin x,\ g(x)=\cos x, \ F(x)=\frac{x^2-1}{x^2}  , \ G(x)=\frac{x^2}{x^2-1}  $. Then $F(f)=G(g)=-\tan^2$, while $f,\ g\in {\mathcal A}(\mathbb{C})$. According to the notation of Theorem 1, we have $p=q=2$, hence $k$ must be at most $1$. Actually this is right because $F'(x)=\frac{2}{x^3}$ has no zero.   
\medskip\\{\bf Acknowledgement.}\quad The authors are very grateful to C.-C. Yang for letting us  know the problem of the
decomposition and the first results obtained in this domain by himself and his colleagues P. Li and  H. K. Ha.
\section{The proofs}
In order to prove the theorems we first need some notation and lemmas
\medskip\\{\bf Notation.} \quad 
Let $\phi , \psi $ be functions from $\mathbb{R}_+$ to $\mathbb{R}_+$. We put $\phi \  \widetilde <   \psi$ if there exists
another  function
$\theta$ from $\mathbb{R}_+$ to $\mathbb{R}_+$ and a subset $H$ of $\mathbb{R}_+$ of zero Lebesgue measure such that $\lim_{r\rightarrow\infty,\;r\notin H}\frac{\theta(r)}{\psi(r)}=0$ and such that $\phi(r)\leq \psi(r)+\theta(r)$ (we notice that
$\widetilde <$ is a preorder). \medskip\\
In the same way, we put  $\phi \  \widetilde -   \psi$ if there exists
another  function $\theta$ from $\mathbb{R}_+$ to $\mathbb{R}_+$ and a subset $H$ of $\mathbb{R}_+$ of zero Lebesgue measure such that $\lim_{r\rightarrow \infty,\;r\notin H}\frac{\theta(r)}{\psi(r)}=0$  and such that $\phi(r)=  \psi(r)+\theta(r)$ (we notice that $\widetilde
-$ is an equivalence relation). 
\medskip\\ Lemma 1 is shown in \cite{6} in a p-adic field of any characteristic $p\geq 0$. Actually it holds in any algebraically closed field:
\medskip \\{\bf Lemma 1.}\quad {\it Let $\mathbb{E}$ be an algebraically closed   field, let $F=\frac{A}{B} , \ G=\frac{C}{D} \in \mathbb{E}(x)$ with $gcd(A,B)=gcd(C,D)=1$,  satisfy Condition $M$ for  $c_1,\dots, c_k$    and let $\deg(G)=q$. Then for each
$j=1,\dots,k$, we have factorizations  of the form $F(x)-F(c_j)= (x-c_j)^{s_j}R_j(x)$ with $s_j\geq 2,\ R_j(c_j)\neq 0$  and $(G(x)-F(c_j))D(x)=\prod_{l=1}^q(x-b_{j,l})$ where the $b_{j,l}$ are $qk$ distincts elements.}
\medskip\\
Lemma 2 is classical (\cite{4},\cite{5},\cite{7}):
\medskip\\{\bf Lemma 2.}\quad {\it Let $R\in \mathbb{C}(x)$ and let $f\in {\mathcal M}(\mathbb{C})$. Then $T(r,R(f)) \widetilde -\deg(R)T(r,f)$. }
\medskip\\ Lemma 3 is basic:
\medskip\\{\bf Lemma 3.}\quad {\it Let $\mathbb{E}$ be a  field, let $F=\frac{A}{B} \in \mathbb{E}(x)$,  let $c$ be a  zero of
$F'$ such that $F(c)\neq 0$. Then $c$ is a  zero of
$\left( \frac{1}{F}\right)'$, with the same multiplicity.}
\medskip\\{\bf Lemma 4.}\quad {\it Let $\mathbb{E}$ be an algebraically closed   field, let $F=\frac{A}{B} , \ G=\frac{C}{D} \in \mathbb{E}(x)$ with $gcd(A,B)=gcd(C,D)=1$,  satisfy Condition
$M'$ for  $c_j, \ 1\leq j\leq k$. Then $\frac{1}{F},\frac{1}{G}$ also satisfy Condition $M'$ for $c_1,\dots,c_k$.}
\begin{proof} By Lemma 3, $c_1,\dots,c_k$ are zeros of $\left( \frac{1}{F}\right)'$. Next, Conditions 1), 2), 3)
are obviously satisfied by $\frac{1}{F},\frac{1}{G}$. So, we only have to show Condition 4'). Let $u$ be a  zero of 
 $D'-\frac{1}{F(c_j)}C'$. First, $u$ 
is a zero of $C'-F(c_j)D'$, hence by Condition 4') satisfied by $F,\ G$, we have  $F(c_j)\neq G(u)$, therefore $\frac{1}{F(c_j)}\neq \frac{1}{G(u)}$. Moreover, by Condition 4') we have $D(u)C(u)\neq 0$ which is symmetric with regard to $C$ and
$D$.  This ends the proof of Lemma 4.\end{proof}
Let us recall an opportune form of  the  Nevanlinna Second Main Theorem (\cite{3}, \cite{4}, \cite{5}, \cite{7}):
\medskip\\ {\bf Theorem N.}\quad {\it Let $f\in {\mathcal M}(\mathbb{C})$. Let $b_1,\dots,b_n\in \mathbb{C}$. Then} $$(n-1)T(r,f)\widetilde <
\sum_{j=1}^n\overline Z(r,f-b_j)+\overline N(r,f).$$  
\medskip\\ {\it Proof of the Theorems 1, 2 and 3.}\quad     In the three Theorems, we put $F=\frac{A}{B},\ G=\frac{C}{D}$ with$gcd(A,B)=gcd(C,D)=1$  and $C, D$
monic.   By Lemma 1, for each $j=1,\dots,k$ we can factorize $F(x)-F(c_j)$ in the form 
\begin{equation}
(x-c_j)^{s_j}R_j(x)
\end{equation}
with $s_j\geq 2$ and $R_j(c_j)\neq 0$ and similarly 
\begin{equation}
(G(x)-F(c_j))D(x)=\prod_{l=1}^q(x-b_{j,l})
\end{equation}
whereas the  $b_{j,l}$ are $qk$ distinct elements.  Thus we obtain 
\begin{equation}
F(f)-F(c_j)=(f-c_j)^{s_j}R_j(f)=\frac{1}{D(g)}\prod_{l=1}^q(g-b_{j,l})=G(g)-F(c_j)
\end{equation}
Now, by Theorem N  we have 
\begin{equation}
(qk-1)T(r,g)\  \widetilde{< }\sum_{j=1}^k \sum_{l=1}^q \overline Z(r, g-b_{j,l}) +\overline N(r,g)
\end{equation}
Thanks to (2), for each fixed $j$  we obtain
\begin{equation}
\overline Z(r,
(G(g)-F(c_j))D(g))=\sum_{l=1}^q\overline Z(r,g-b_{j,l})
\end{equation}
Then, by inserting $g$ in (4) and in (5) we obtain
\begin{equation}
(qk-1)T(r,g)\widetilde{< }
\sum_{j=1}^k\overline Z(r, (G(g)-F(c_j))D(g))+\overline N(r,g)
\end{equation}
Now, by Lemma 2 we have
$T(r,g)\widetilde  -   \frac{p}{q}T(r,f)$. Let $R_j(x)=\frac{A_j(x)}{B_j(x)}$: of course  we have  $\deg(A_j)\leq p-s_j$.  By applying Lemma 2,  we obtain 
\begin{equation}
\overline Z(r,R_j(f))+\overline N(r,f)=\overline Z(r, A_j(f)) +\overline N(r,f) \widetilde < (p-s_j)T(r,f)+\overline N(r,f)
\end{equation}
Now, assume the hypothesis of Theorem 1.  Since $f$ and $g$ are entire and since  on one hand, $A$ and $B$ have no common zero and on the other hand,  $C$ and $D$ have no common zero, we can see that by  (3), for each $j=1,\dots,k$,
any zero of $D(g)$ is a  pole of $F(f)$ of same order, and therefore is a pole of 
$R_j(f)$ of same order. Therefore $\overline Z(r,R_j(f)D(g))=\overline Z(r,R_j(f))$. Consequently, by (6) and (7), taking into account that $f, \ g$ are entire, we can get
\begin{eqnarray}\nonumber
\overline Z(r, (G(g)-F(c_j))D(g))=\overline Z(r, (f-c_j)^{s_j}R_j(f)D(g))\leq\\\nonumber
\leq \overline Z(r, (f-c_j))+\overline Z(r,
R_j(f))\widetilde < T(r,f)+(p-s_j)T(r,f)
\end{eqnarray}
 Summing up on $j$, by (6) we obtain $ (qk-1)T(r,g)\leq T(r,f)(k(1+p)-\sum_{j=1}^ks_j).$ On the other hand,  by Lemma 2 we have $T(r,f)\widetilde - \frac{q}{p}T(r,g)$, therefore \[(qk-1)pT(r,f)\widetilde < qT(r,f)(kp-\sum_{j=1}^k(s_j-1))\] Since $T(r,f)$ is unbounded, we can derive that $ (qk-1)p \leq   q(kp-\sum_{j=1}^k(s_j-1)) $ hence $\sum_{j=1}^k(s_j-1) \leq p$. But by hypothesis,  $k\leq \sum_{j=1}^k(s_j-1)$, thereby  $kq\leq p$. 
\medskip\\Next, assume the hypothesis of Theorem 2.  We can find a   constant
 $h\in \mathbb{C}$ such that $F(c_j)+h\neq 0\ \forall j=1,\dots,k$. Let $\Phi(x)=F(x)+h,\ \Psi(x)=G(x)+h$.
Thus we have $\Phi(f)=\Psi(g)$, and $\Phi , \ \Psi$ satisfy the hypothesis of Theorem 2, in place of $F,\ G$, respectively. Consequently, we can assume that $F(c_j)\neq 0\ \forall j=1,\dots,k$, without loss of generality. Now, we  have  \[\overline Z(r,R_j(f)D(g))\leq \overline Z(r,R_j(f))+\overline Z(r,D(g))\] Consequently, by (6) and (7),  we can get
\begin{eqnarray}
\overline Z(r,
(G(g)-F(c_j))D(g))=\overline Z(r, (f-c_j)^{s_j}R_j(f)D(g))\leq\\\nonumber\leq\overline Z(r, (f-c_j))+\overline Z(r,  R_j(f))+\overline Z(r, D(g))
\end{eqnarray}
Let $t=\gamma(D)$ and let $\delta_u \ (1\leq u\leq t)$ be the distinct zeros of $D$. 
We notice that $\overline Z(r, D(g))\leq \sum_{u=1}^t\overline Z(r, g-\delta_u)\leq tT(r,g)$. Consequently, by (8) we obtain  
\[\overline Z(r,
(G(g)-F(c_j))D(g)) \widetilde < T(r,f)+(p-s_j)T(r,f)+\gamma (D)T(r,g)\widetilde <
(p-(s_j-1)+\frac{p\gamma(D)}{q})T(r,f)\]
\medskip\\Summing up on
$j$, by (6) we obtain \[(qk-1)T(r,g)\leq T(r,f)(k(1+p)-\sum_{j=1}^ks_j)\] therefore 
\[(qk-1)pT(r,f)\widetilde < qT(r,f)(kp-\sum_{j=1}^k(s_j-1)+\frac{kp\gamma(D)}{q})\] Since $T(r,f)$ is unbounded, we can derive that $ (qk-1)p \leq   q(kp-\sum_{j=1}^k(s_j-1))+kp\gamma(D) $ therefore   $kq\leq
q\sum_{j=1}^k(s_j-1)\leq  p(k\gamma(D)+1)$. \\ Now, assume the hypothesis of Theorem 3.  Thanks to Condition M', by Lemma 4   we can replace $F,\  G$ by $\frac{1}{F} $ and $\frac{1}{G}$ respectively,   because they also satisfy Condition M',   hence  we obtain  $kq\leq  p(k\gamma(C)+1)$,
consequently $kq\leq  p(k\lambda(G)+1)$.\ep

\end{document}